\documentclass{article}

\usepackage{latexsym}
\usepackage{amsmath}
\usepackage{amssymb}

\title{A Note on Walk versus Wait: Lazy Mathematician Wins\footnote{J.G. Chen, S.D. Kominers, and R.W. Sinnott. Walk versus wait: The lazy mathematician wins. arXiv.org Mathematics, January 2008. http://arVix.org/abs/0801.0297}}
\author{Ramnik Arora\\Indian Institute of Technology Kanpur}
\date{\today}

\begin{document}
\maketitle
It seems that the distance term in the given equation $4$ and the expression above it in the article \cite{WvW} is not generalised enough and is partially incorrect:
The expression and the equation are reproduced here:
$$
\frac{d_2}{v_w} + \underbrace{\int^{t_w}_{0}( \frac{1}{t_b} ( \frac{ d - d_2 }{v_b}))dt}_1 + \underbrace{ ( 1 - \int^{t_w}_{0}p(t)dt)(\frac{d}{v_w} + t_w ) }_2
$$
and
$$
\underbrace{ \int^{t_w}_{0}( \frac{1}{t_b} ( \frac{ d - d_2 }{v_b} + t))dt }_1 + \underbrace{ ( 1 - \int^{t_w}_{0}p(t)dt)(\frac{d}{v_w} + t_w ) }_2 = \frac{d - d_2}{v_w}\\
$$

\section*{Distance}
The distance ($d$) and the time used in the term 2 is not correct. Thus, on correcting the distance, the $2^{nd}$ term would be:

$$ 
(1 - \int^{t_w}_{0}p(t)dt)(\frac{d - d_2}{v_w} + t_w ) 
$$

It is because if we use $d$ instead of $d-d_2$, then we are in effect double counting the distance $d_2$
when the mathematician is walking all the way to the destination.

\section*{Time}
The probability density function of the bus reaching the second bus stop will not be only $p(t)$ as is suggested in the
article \cite{WvW}. It can be seen that the probability of the bus reaching the second bus stop after waiting time $t=t_\circ$ is actually
given by $p( t_\circ - \frac{d_2}{v_b} + \frac{d_2}{v_w} )$.

This can be explained by seeing that we have walked for $\frac{d_2}{v_w}$ hours before coming to the second bus stop and
that for a bus to reach the second bus stop at time $t_\circ'$, it needs to be at the first bus stop at time $t_\circ' - \frac{d_2}{v_b}$. 
This is keeping in mind that $p(t)$ was defined as the probability of the bus arriving at the first bus stop at time $t$.

\section*{Generalisation}

Also, it seems that in equation 4 (the second equation here), the result can be generalised by using $p(t_{corrected})$ probability distribution
instead of the very specific $\frac{1}{t_d}$.
This makes the $1^{st}$ term:
$$
\int^{t_w}_{0} p(t_{corrected})(\frac{d - d_2}{v_b} + t)dt
$$

where $t_{corrected}$ is defined as $t - \frac{d_2}{v_b} + \frac{d_2}{v_w}$, as is suggested above.

\section*{Residual term}
Moreover, it seems, there will be another factor that has not been considered in the expressions/equations.
Appearance of the bus at stop $1$ is a sort of a periodic event (stochastic process).
Hence, in the uniform distribution case, there will be significant dependence of the waiting time on
the moment the bus passes us by while we are en route to the second bus stop. This term would come out
to be (under the assumption of uniform distribution):

$$
\int^{\frac{d_2}{v_w}}_{0}\frac{1}{t_b}[ \underbrace{( t_b - t )}_1 - \underbrace{( \frac{d_2 - v_wt}{v_w} )}_2 ]dt
$$

This is the expected time that he has to wait additionally if the bus happens to pass him at while he was on his way.
In case of uniform distribution we know that one bus is expected from $[0,t_b]$ and the next from $[t_b,2t_b]$.
The $1^{st}$ term is the 'dead' time for which the next bus is not expected and the $2^{nd}$ term is the time that he
spends walking to the bus stop, which is subtracted from the $1^{st}$ term.

This is under the assumption that $\frac{d_2}{v_w} < t_b$ as if this is not the case, then the Mathematician would always
choose to wait for the bus as it will necessarily pass him before he gets to his destination.

This is just for the simple case of uniform distribution, and the result for the general distribution can also be worked out
using similar arguments.\\

Nevertheless, these changes will not make a difference to the validity of the result. 

\section*{Acknowledgement}
I have to thank Utkarsh Upadhyay for his helpful comments and backspaces.

\end{document}